\documentclass[final,leqno]{article}
\usepackage{epsfig,psfrag,amsmath,amssymb,url}
\usepackage[letterpaper=true,colorlinks=true,linkcolor=red,citecolor=red]{hyperref}

\newtheorem{mytheorem}{Theorem}
\newtheorem{mylemma}{Lemma}

\newtheorem{myremark}[mytheorem]{Remark}

\newcommand{\N}{{\sf N\hspace*{-1.0ex}\rule{0.15ex}{1.3ex}\hspace*{1.0ex}}}

\newcommand{\BA}{\begin{eqnarray}}
\newcommand{\EA}{\end{eqnarray}}
\newcommand{\BE}{\begin{equation}}
\newcommand{\EE}{\end{equation}}
\newcommand{\ba}{\begin{array}}
\newcommand{\ea}{\end{array}}
\newcommand{\baa}{\begin{eqnarray*}}
\newcommand{\eaa}{\end{eqnarray*}}
\newcommand{\be}{\begin{equation}}
\newcommand{\ee}{\end{equation}}

\def\F{{\cal F}}

\def\N{{\cal N}}

\def\u{{\bf u}}

\def\Rtilde{{\widetilde{R}}}
\def\Btilde{{\widetilde{B}}}

\def\Atilde{{\widetilde{A}}}
\def\Wtilde{{\widetilde{W}}}
\def\Vtilde{{\widetilde{V}}}

\def\>{\raisebox{-1ex}{$\; \stackrel{\textstyle >}{\sim } \; $}}
\def\<{\raisebox{-1ex}{$ \; \stackrel{\textstyle <} {\sim } \; $}}
\def\myalpha{{\alpha }}
\newcommand{\vvec}[1]{{\mathbf{#1}}}

\newcommand{\EQ}[1]{(\ref{equation:#1})}
\newcommand{\LA}[1]{\ref{lemma:#1}}

\def\beginproof{\indent {\it Proof:~}}
\def\endproof{$\Box$}

\title{A FETI-DP type domain decomposition algorithm for
 three-dimensional  incompressible Stokes equations}

\author{Xuemin Tu\thanks{Department of Mathematics, University of Kansas, 1460 Jayhawk Blvd, Lawrence, KS 66045-7594,  {\tt xtu@math.ku.edu}, {\tt http://www.math.ku.edu/$\sim$xtu/}. This author's work was supported in part by National Science Foundation contract DMS-1115759.} \and Jing Li\thanks{Department of Mathematical Sciences, Kent State University, Kent, OH 44242, {\tt li@math.kent.edu}, {\tt http://www.math.kent.edu/$\sim$li/}.} }

\begin{document}

\maketitle
\pagestyle{myheadings} \thispagestyle{plain} \markboth{XUEMIN TU AND
  JING LI}{DOMAIN DECOMPOSITION FOR INCOMPRESSIBLE STOKES}

\begin{abstract}
The FETI-DP algorithms, proposed by the authors in [{\em SIAM J. Numer. Anal.}, 51 (2013), pp.~1235--1253] and [{\em Internat. J. Numer. Methods Engrg.}, 94 (2013), pp.~128--149] for
solving incompressible Stokes equations, are extended to
three-dimensional problems. A new analysis of the condition number
bound for using the Dirichlet preconditioner is given. An advantage of this new analysis is that the numerous coarse level velocity components, required in the previous analysis to enforce the divergence free subdomain boundary velocity conditions, are no longer needed. This greatly reduces the size of the coarse level problem in the algorithm, especially for three-dimensional problems. The coarse level velocity space can be chosen as simple as for solving scalar elliptic problems corresponding to each velocity component. Both Dirichlet and lumped preconditioners are analyzed using a same framework in this new analysis. Their condition number bounds are proved to be independent of the number of subdomains for fixed subdomain problem size. Numerical experiments in both two and three dimensions demonstrate the convergence rate of the algorithms.

\end{abstract}

%\begin{keywords}
{\bf Keywords}
domain decomposition, incompressible Stokes, FETI-DP, BDDC, divergence free

%\end{keywords}

%\begin{AMS}
{\bf AMS} 
65F10, 65N30, 65N55
%\end{AMS}

\section{Introduction}
Mixed finite elements are often used to solve incompressible Stokes
and Navier-Stokes equations. Continuous pressures have been used in
many mixed finite elements, e.g., the well known Taylor-Hood finite
elements \cite{Taylor}.  However, most domain decomposition methods require that the pressure be discontinuous, when they are used to solve the indefinite linear systems arising from such mixed finite element discretizations;  see, e.g., \cite{Doh04, Doh09, Doh10, paulo2003, kim06, Kla98, li05, li06, 8pavwid, LucaOlofStef, Tu:2005:BPP, Tu:2005:BPD}.
Several domain decomposition algorithms allow  to use continuous pressures, e.g.,  Klawonn and Pavarino~\cite{Kla98}, Goldfeld~\cite{pauloPhD},  \v{S}\'istek
{\em et. al.}~\cite{sis11},  Benhassine and Bendali~\cite{ben10}, and Kim and Lee~\cite{kim12}. But the convergence rate analysis of those approaches
cannot be applied to the continuous pressure case due to the indefiniteness of the linear systems;
such difficulty can often be removed conveniently when discontinuous
pressures are used in the discretization.
% for domain decomposition methods.

Recently, the authors  \cite{li12, Tu:2012:Stokes} proposed and
analyzed a FETI-DP (Dual-Priaml Finite Element Tearing and
Interconnecting method) type domain decomposition algorithm for
solving  the incompressible Stokes equation in two
dimensions. Both discontinuous and continuous pressures can be used in the
mixed finite element discretization. In both cases, the indefinite
system of linear equations can be reduced to a symmetric positive
semi-definite system. Therefore, the preconditioned conjugate gradient method
can be applied and a scalable convergence rate of the algorithm has been proved.

The lumped and Dirichlet preconditioners have been studied in \cite{li12} and \cite{Tu:2012:Stokes}, respectively.
%Compared with the Dirichlet preconditioner, the lumped preconditioner is less effective in the reduction of the iteration count, but it is also less expensive. The lumped preconditioner requires only subdomain matrix and vector products, while the  Dirichlet preconditioner requires solving subdomain systems of equations.
For the lumped preconditioner it was shown both experimentally and analytically in \cite{li12}, that the coarse level space can be chosen the same as for solving scalar elliptic problems corresponding to each velocity component to achieve a scalable convergence rate. Similar observations for the lumped preconditioner have also been pointed out earlier by Kim and Lee~\cite[with Park]{kim11, kim102, kim10}, even though their studies are only for using discontinuous pressures.

For the Dirichlet preconditioner studied in \cite{Tu:2012:Stokes}, a distinctive feature is the application of subdomain discrete harmonic extensions in the preconditioner.  In other existing FETI-DP and BDDC (Balancing Domain Decomposition by Constraints) algorithms, cf. \cite{li05, li06}, subdomain discrete Stokes extensions have been used and the coarse level velocity space has to contain sufficient components to enforce divergence free subdomain boundary velocity conditions.
Those complicated and numerous coarse level velocity components, especially for three-dimensional problems as discussed in \cite{li06}, are not needed for the implementation of the Dirichlet preconditioner in \cite{Tu:2012:Stokes}. But they are still required in \cite{Tu:2012:Stokes} just for the analysis, where subdomain Stokes extensions were used, to obtain a scalable condition number bound.

In this paper,  we provide a new analysis for the algorithms in
\cite{li12, Tu:2012:Stokes}, which can analyze both lumped and Dirichlet preconditioners in a same framework. It does not use any subdomain Stokes extensions and those additional coarse level velocity
components to enforce divergence free subdomain boundary velocity conditions are no longer needed. For both lumped and Dirichlet preconditioners, the coarse level space can be chosen as simple as for solving scalar elliptic problems corresponding to each velocity component.
This  greatly simplifies the requirements on the coarse level space for the case of Dirichlet preconditioner, especially in  three dimensions.  This paper is presented in the context of solving three-dimensional problems; the same approach can be applied to two-dimensional problems as well.

The remainder of this paper is organized as follows. The finite
element discretization of the incompressible Stokes equation is
introduced in Section \ref{section:FEM}. A domain decomposition
approach is described in Section~\ref{section:DDM}, and the system is reduced to a symmetric positive semi-definite problem in Section~\ref{section:Gmatrix}. A few preliminary results used in the
condition number bound estimates are given in Section
\ref{section:techniques}. The lumped and Dirichlet preconditioners are introduced in Section
\ref{section:jump}, and the condition number bounds of the
preconditioned systems are established in Section~\ref{section:convergence}. At
the end, numerical results of solving the incompressible Stokes equation in both two and three dimensions are given in Section~\ref{section:numerics} to demonstrate the convergence rate of the algorithm.

\section{Finite element discretization}
\label{section:FEM}

We consider solving the following incompressible Stokes problem on a
bounded, three-dimensional polyhedral domain $\Omega$  with a
Dirichlet boundary condition,
\begin{equation}
\label{equation:Stokes}
\left\{
\begin{array}{rcll}
-\Delta {\bf u}^* + \nabla p^* & = & {\bf f}, & \mbox{ in } \Omega \mbox{ , } \\
-\nabla \cdot {\bf u}^*        & = & 0, & \mbox{ in } \Omega \mbox{ , } \\
{\bf u}^*                      & = & {\bf u}^*_{\partial \Omega}, & \mbox{ on } \partial \Omega \mbox{ , }\\
\end{array}\right.
\end{equation}
where the boundary velocity ${\bf u}^*_{\partial \Omega}$ satisfies the compatibility
condition $\int_{\partial \Omega} {\bf u}^*_{\partial \Omega} \cdot {\bf n} = 0$. For
simplicity, we assume that ${\bf u}^*_{\partial \Omega} = {\bf 0}$ without losing any
generality.

The weak solution of \EQ{Stokes} is given by: find $\vvec{u}^* \in
\left(H^1_0(\Omega)\right)^3 = \{ \vvec{v} \in (H^1(\Omega))^3 ~
\big| ~ \vvec{v} = \vvec{0} \mbox{ on }
\partial \Omega \}$ and $p^* \in
L^2(\Omega)$, such that
\begin{equation}
\label{equation:bilinear} \left\{
\begin{array}{lcll}
a(\vvec{u}^*, \vvec{v}) + b(\vvec{v}, p^*) & = & (\vvec{f}, \vvec{v}),
& \forall \vvec{v}\in \left(H^1_0(\Omega)\right)^3 , \\ [0.5ex]
b(\vvec{u}^*, q) & = & 0, & \forall q \in L^2(\Omega) \mbox{ , }
\\
\end{array} \right.
\end{equation}
where
$
a(\vvec{u}^*, \vvec{v})= \int_{\Omega} \nabla{\bf u}^* \cdot \nabla{\bf v}, ~
b(\vvec{u}^*,q) = -\int_{\Omega} (\nabla \cdot \vvec{u}^*) q, ~
(\vvec{f}, \vvec{v}) = \int_{\Omega} \vvec{f} \cdot \vvec{v}.
$
We note that the solution of \EQ{bilinear} is not unique, with the pressure $p^*$ different up to an additive constant.

A mixed finite element is used to solve \EQ{bilinear}. In this paper we apply a mixed finite element with continuous pressures, e.g., the Taylor-Hood type mixed finite elements. The same algorithm and analysis can be applied to mixed finite elements with discontinuous pressures as well; see \cite{Tu:2012:Stokes}. Denote the velocity finite element space by $\vvec{W} \subset \left(H^1_0(\Omega)\right)^3$, and the
pressure finite element space by $Q \subset L^2(\Omega)$.
The finite element solution $(\vvec{u}, p) \in \vvec{W} \bigoplus Q$ of \EQ{bilinear} satisfies
\begin{equation}
\label{equation:matrix} \left[
\begin{array}{cccc}
A     &  B^T \\
B     &  0   \\
\end{array}
\right] \left[
\begin{array}{c}
{\bf u}     \\
p   \\
\end{array}
\right] = \left[
\begin{array}{l}
{\bf f}        \\
0      \\
\end{array}
\right] ,
\end{equation}
where $A$, $B$, and $\vvec{f}$  represent respectively the restrictions of $a( \cdot , \cdot )$, $b(\cdot, \cdot )$ and
$(\vvec{f} , \cdot)$ to the finite-dimensional spaces $\vvec{W}$ and $Q$. We use the same notation in this paper to represent both a finite element function and the vector of its nodal values.

The coefficient matrix in \EQ{matrix} is rank deficient even though
$A$ is symmetric positive definite. $Ker(B^T)$, the kernel of $B^T$,
contains all constant pressures in $Q$.  $Im(B)$, the range of
$B$,  is orthogonal to $Ker(B^T)$ and consists of all vectors in $Q$ with zero average.  For a general right-hand side vector $({\bf f}, ~ g)$ in \EQ{matrix}, the existence of solution requires that $g \in Im(B)$, i.e., $g$ has zero average; for the right-hand side given in \EQ{matrix}, $g = 0$ and the solution always exists. When the pressure is considered in the quotient space $Q/Ker(B^T)$, the solution is unique.  In this paper, when $q \in Q/Ker(B^T)$, we always assume that $q$ has zero average.

Let $h$ represent the characteristic diameter of the mixed elements. We assume that the mixed finite element space $\vvec{W} \times Q$, is inf-sup stable in the sense that there exists
a positive constant $\beta$, independent of $h$, such that
\begin{equation}
\label{equation:infsupMatrix} \sup_{{\bf w} \in {\bf W}}
\frac{\left< q, B \vvec{w} \right>^2}{\left< \vvec{w},  A \vvec{w} \right>} \geq \beta^2 \left< q, Z q \right>, \hspace{0.5cm} \forall q \in Q/Ker(B^T),
\end{equation}
cf.~\cite[Chapter III, \S 7]{braess}.  Here, as always used in this paper, $\left< \cdot, \cdot \right>$
represents the inner (or semi-inner) product of two vectors. The matrix $Z$ represents the mass matrix defined on the pressure
finite element space $Q$, i.e., for any $q \in Q$, $\|q\|_{L^2}^2 =
\left< q, Z q \right>$. It is easy to see, cf.~\cite[Lemma
B.31]{Toselli:2004:DDM}, that $Z$ is spectrally equivalent to $h^3 I$
for three-dimensional problems, i.e., there exist positive constants $c$ and $C$, such that
\be
\label{equation:massmatrix}
c h^3 I \leq Z \leq C h^3 I,
\ee
where $I$ represents the identity matrix. Here, as in other places of this paper, $c$ and $C$  represent generic positive constants which are independent of $h$ and the subdomain diameter $H$ (described in the following section).

\section{A non-overlapping domain decomposition approach}
\label{section:DDM}

The domain $\Omega$ is decomposed into $N$ non-overlapping polyhedral
subdomains $\Omega_i$, $i = 1, 2, ..., N$. Each subdomain is the union
of a bounded number of elements, with the diameter of the subdomain in
the order of $H$. We use $\Gamma$ to represent the subdomain interface
which contains all the subdomain boundary nodes shared by neighboring
subdomains; we assume that the subdomain meshes have matching nodes
across $\Gamma$.
$\Gamma$ is composed of subdomain faces, which are regarded as open subsets of
$\Gamma$ shared by two subdomains, subdomain edges, which are regarded
as open subsets of $\Gamma$ shared by more than two subdomains, and of the subdomain vertices, which are end points of edges.

The velocity and pressure finite element spaces ${\bf W}$ and $Q$ are decomposed into
\[
{\bf W} = {\bf W}_I \bigoplus {\bf W}_{\Gamma}, \quad
Q = Q_I \bigoplus Q_\Gamma,
\]
where ${\bf W}_I$ and $Q_I$ are direct sums of independent subdomain interior velocity spaces ${\bf W}^{(i)}_I$, and interior pressure spaces $Q^{(i)}_I$, respectively, i.e.,
$$
{\bf W}_I = \bigoplus_{i=1}^{N}{\bf W}^{(i)}_I, \quad Q_I =
\bigoplus_{i=1}^{N}Q^{(i)}_I.
$$
${\bf W}_{\Gamma}$ and $Q_\Gamma$ are subdomain interface velocity and pressure spaces, respectively. All
functions in ${\bf W}_{\Gamma}$  and $Q_\Gamma$ are continuous across $\Gamma$; their degrees of freedom are shared by neighboring subdomains.

To formulate the domain decomposition algorithm, we introduce a partially sub-assembled subdomain interface velocity space $\vvec{\Wtilde}_{\Gamma}$,
\[
\vvec{\Wtilde}_{\Gamma} = \vvec{W}_{\Delta} \bigoplus \vvec{W}_{\Pi}
= \left( \bigoplus_{i=1}^N \vvec{W}^{(i)}_\Delta \right) \bigoplus \vvec{W}_{\Pi}.
\]
$\vvec{W}_{\Pi}$ is the continuous, coarse level, primal velocity space which is typically spanned by subdomain vertex nodal basis functions, and/or by interface edge/face-cutoff functions with constant nodal values on each edge/face, or with values of positive weights on these edges/faces. The primal, coarse level velocity degrees of freedom are shared by neighboring subdomains. The complimentary space $\vvec{W}_{\Delta}$ is the direct sum of independent subdomain dual interface velocity spaces $\vvec{W}_{\Delta}^{(i)}$, which correspond to the remaining subdomain interface velocity degrees of freedom and are spanned by basis functions which vanish at the primal degrees of freedom. Thus, an element in $\vvec{\Wtilde}_{\Gamma}$ typically has a continuous primal velocity component and a discontinuous dual velocity component.

It is well known that, for domain decomposition algorithms, the coarse space $\vvec{W}_{\Pi}$ should be sufficiently rich to achieve a scalable convergence rate. On the other hand, a large coarse level problem will certainly degrade the parallel performance of the algorithm. Therefore it is important to keep the size of the coarse level problem as small as possible.
When the Dirichlet preconditioner was used in the FETI-DP algorithm for solving incompressible Stokes equations \cite{li05} and similarly in the BDDC algorithm \cite{li06}, subdomain discrete Stokes extensions were used and $\vvec{W}_{\Pi}$ has to contain sufficient subdomain interface components such that functions in $\vvec{W}_{\Delta}$ have zero flux across the subdomain boundaries.
Such requirements lead to a large coarse level velocity space, especially for three-dimensional problems, cf. \cite{li06}.

In \cite{Tu:2012:Stokes}, a FETI-DP type algorithm
is proposed for solving two-dimensional incompressible Stokes
problems. A distinctive feature of the
Dirichlet preconditioner used in that algorithm is the application of
subdomain discrete harmonic extensions, instead of subdomain discrete
Stokes extensions. As a result, the divergence free subdomain boundary
velocity conditions are not needed in that algorithm.  However, the
analysis, given in \cite{Tu:2012:Stokes} for the Dirichlet
preconditioner, still uses subdomain Stokes extensions and requires the
same type coarse level velocity space as discussed in \cite{li06} to establish  a scalable condition number bound estimate. In this paper, a new analysis is offered and it is sufficient for $\vvec{W}_{\Pi}$ to be spanned just by the subdomain vertex nodal basis functions and subdomain edge-cutoff functions corresponding to each velocity component, as for solving three-dimensional scalar elliptic problems, cf.~\cite[Section~6.4.2]{Toselli:2004:DDM}.

The functions ${\bf w}_{\Delta}$ in ${\bf W}_{\Delta}$ are in general
not continuous across $\Gamma$. To enforce their continuity, we define
a Boolean matrix $B_\Delta$ of the form
\[
B_\Delta = \left[ B_\Delta^{(1)} ~~~ B_\Delta^{(2)} ~~~ \cdots ~~~ B_\Delta^{(N)} \right],
\]
constructed from $\{0,1,-1\}$. On each
row of $B_\Delta$, there are only two nonzero entries, $1$ and $-1$,
corresponding to one velocity degree of freedom shared by two neighboring subdomains, such that
for any ${\bf w}_{\Delta}$ in ${\bf W}_{\Delta}$, each row of
$B_\Delta {\bf w}_{\Delta} = 0$ implies that these two degrees of
freedom from the two neighboring subdomains be the same. We note that,
in three dimensions, a velocity degree of freedom on a subdomain edge
is shared by more than two subdomains, e.g., by four subdomains. In this case,  a minimum of three continuity constraints can be applied to enforce the continuity of this velocity degree of freedom among the four subdomains, which corresponds to the use of non-redundant Lagrange multipliers. In this paper, the fully redundant Lagrange multipliers are used, which means, e.g., for a subdomain edge velocity degree of freedom shared by four subdomains, six Lagrange multipliers are used to enforce all the six possible continuity constraints among them, cf.~\cite[Section~6.3.1]{Toselli:2004:DDM}.

We denote the range of $B_\Delta$ applied on ${\bf W}_{\Delta}$ by $\Lambda$, the vector space of the Lagrange multipliers. Solving the original fully assembled linear system~\EQ{matrix} is then equivalent to: find
$\left( {\bf u}_I, ~p_I, ~{\bf u}_{\Delta}, ~{\bf u}_{\Pi}, ~p_{\Gamma}, ~\lambda \right) \in
{\bf W}_I \bigoplus Q_I \bigoplus {\bf W}_{\Delta} \bigoplus {\bf W}_\Pi \bigoplus Q_\Gamma \bigoplus \Lambda$, such that
\be
\label{equation:bigeq}
\left[
\begin{array}{cccccc}
A_{II}      & B_{II}^T       & A_{I \Delta}       & A_{I \Pi}       & B_{\Gamma I}^T     &  0           \\[0.8ex]
B_{II}      & 0              & B_{I \Delta}       & B_{I \Pi}       & 0                  &  0           \\[0.8ex]
A_{\Delta I}& B_{I \Delta} ^T& A_{\Delta\Delta}   & A_{\Delta \Pi}  & B_{\Gamma \Delta}^T&  B_{\Delta}^T\\[0.8ex]
A_{\Pi I}   & B_{I \Pi}^T    & A_{\Pi \Delta}     & A_{\Pi \Pi}     & B_{\Gamma \Pi}^T   &  0           \\[0.8ex]
B_{\Gamma I}& 0              & B_{\Gamma \Delta}  & B_{\Gamma \Pi}  & 0                  &  0           \\[0.8ex]
0           & 0              & B_{\Delta}         & 0               & 0                  &  0
\end{array}
\right]
\left[ \begin{array}{c}
{\bf u}_I        \\[0.8ex]
p_I              \\[0.8ex]
{\bf u}_{\Delta} \\[0.8ex]
{\bf u}_{\Pi}    \\[0.8ex]
p_{\Gamma}     \\[0.8ex]
\lambda
\end{array} \right] =
\left[ \begin{array}{l}
{\bf f}_I        \\[0.8ex]
0                \\[0.8ex]
{\bf f}_{\Delta} \\[0.8ex]
{\bf f}_\Pi      \\[0.8ex]
0                \\[0.8ex]
0
\end{array} \right] \mbox{ ,  }
\ee
where the sub-blocks in the coefficient matrix represent the restrictions of $A$ and $B$ in~\EQ{matrix} to appropriate subspaces. The leading three-by-three block can be made block diagonal with each diagonal block corresponding to one subdomain.

The coefficient matrix in \EQ{bigeq} is singular. The trivial null space vectors are those with $\lambda$ in the null space of $B_\Delta^T$ and other components zero. Such singularity, due to the rank deficiency of $B_\Delta$, needs not to be worried, since the Lagrange multiplier vector $\lambda$ will be confined in  $\Lambda$, the range of $B_\Delta$.  The only meaningful  basis vector in the null space of \EQ{bigeq} corresponds to the one-dimensional null space of the original incompressible Stokes system~\EQ{matrix}, and is specified in the following lemmas.

We first need to introduce a positive scaling factor $\delta^{\dagger}(x)$ for each node $x$ on $\Gamma$. Let $\N_x$ be the number of subdomains sharing $x$, and we define $\delta^{\dagger}(x) = 1/\N_x$. Given such scaling factors at the subdomain interface nodes, we can define a scaled operator $B_{\Delta, D}$. We note that each row of $B_\Delta$ has only two nonzero entries, $1$ and $-1$, connecting two neighboring subdomains sharing a node $x$ on $\Gamma$. Multiplying each entry by the scaling factor $\delta^{\dagger}(x)$ gives us $B_{\Delta, D}$. Namely
\[
B_{\Delta, D} = \left[ D_\Delta B_\Delta^{(1)} ~~~ D_\Delta B_\Delta^{(2)} ~~~ \cdots ~~~ D_\Delta B_\Delta^{(N)} \right],
\]
where $D_\Delta$ is a diagonal matrix and contains $\delta^{\dagger}(x)$ on its diagonal. We also see from the definition of $B_{\Delta, D}$ that the scalings on all the Lagrange multipliers related to the same subdomain interface node are the same, from which we have the following lemma.

\begin{mylemma}
\label{lemma:BdeltaNull} The null of $B_\Delta^T$ is the same as the null of $B_{\Delta, D}^T$;
the range of $B_\Delta$ is the same as the range of $B_{\Delta, D}$.
\end{mylemma}

The following lemma can be found at \cite[Page 175]{Toselli:2004:DDM}.

\begin{mylemma}
\label{lemma:Bdelta} For any $\lambda \in \Lambda$,
$B_\Delta B_{\Delta, D}^T \lambda =  B_{\Delta, D} B_\Delta^T \lambda = \lambda.$
\end{mylemma}

\begin{mylemma}
\label{lemma:faceIntegral} Let $1_{p_I} \in Q_I$, $1_{p_\Gamma} \in Q_\Gamma$ represent vectors with value $1$ on each entry. Then
\be
\label{equation:aboutlambda}
[B_{I\Delta}^T ~~ B_{\Gamma\Delta}^T]\left[\begin{array}{c}1_{p_I}\\
1_{p_\Gamma}\end{array}\right] = B_\Delta^T \lambda,
\ee
where
\be
\label{equation:lambda}
\lambda = B_{\Delta,D}[B_{I\Delta}^T ~~ B_{\Gamma\Delta}^T]\left[\begin{array}{c}1_{p_I}\\
   1_{p_\Gamma}\end{array}\right] \in \Lambda.
\ee
\end{mylemma}

\beginproof The left side of \EQ{aboutlambda} contains face integrals of the
normal component of the dual subdomain interface velocity finite
element basis functions across the subdomain interface. For a face velocity degree of freedom, which is shared by two
neighboring subdomains, the face integrals of their normal components
on the two neighboring subdomains are negative of each other, since
their normal directions are opposite. This pair of opposite values
can then be represented by the product of $B_\Delta^T$ and a Lagrange
multiplier with value equal to the face integral of the corresponding
basis function.

\begin{figure}[h]
\begin{center}
\input{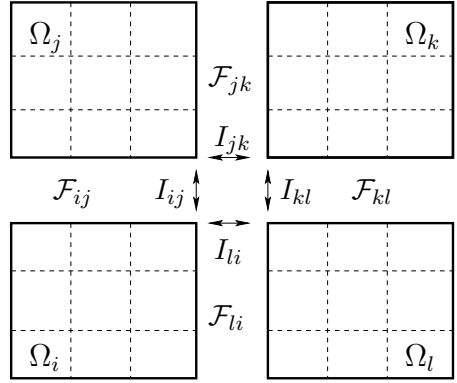}
\end{center}
\caption{Illustration on a subdomain edge interface degree of freedom.}
\label{fig:EdgeNode}
\end{figure}

Now we consider a subdomain edge velocity degree of freedom, which is
shared by more than two subdomains, e.g., by four subdomains
$\Omega_i$, $\Omega_j$, $\Omega_k$, and $\Omega_l$. A two-dimensional
illustration of such an edge node is shown in Figure
\ref{fig:EdgeNode}, where the edge shared by the four subdomains
points outward directly. Denote the four faces having this edge in
common by $\F_{ij}$,  $\F_{jk}$, $\F_{kl}$, $\F_{li}$,  where, e.g.,
$\F_{ij}$ represents the face shared by $\Omega_i$ and $\Omega_j$,
while $\Omega_i$ and $\Omega_k$ have no common face. Denote the
integration of the normal component of this velocity basis function on
these four faces by $I_{ij}$,  $I_{jk}$, $I_{kl}$, $I_{li}$, with a
chosen normal direction for each face, e.g., upward on $\F_{ij}$ and
$\F_{kl}$, to the right on $\F_{jk}$ and $\F_{li}$.  Then the entries
of the left side vector in \EQ{aboutlambda} corresponding to this edge velocity degree of
freedom on the four subdomains $\Omega_i$, $\Omega_j$, $\Omega_k$, and
$\Omega_l$, are $I_{ij} + I_{li}$, $-I_{ij} + I_{jk}$, $-I_{jk} -
I_{kl}$, and $I_{kl} - I_{li}$, respectively. Here two neighboring
subdomains sharing a common face have opposite face integral values on
that face because their normal directions are opposite of each
other. Take $I_{ij}$,  $I_{jk}$, $I_{kl}$, $I_{li}$ as the four
Lagrange multiplier values as illustrated in Figure
\ref{fig:EdgeNode}. Then the four subdomain face integral values
$I_{ij} + I_{li}$, $-I_{ij} + I_{jk}$, $-I_{jk} - I_{kl}$, and $I_{kl}
- I_{li}$, can be represented as the product of corresponding $B_\Delta^T$ with a Lagrange multiplier vector containing these four Lagrange multiplier values and zero elsewhere.

The above has just shown that the left side of \EQ{aboutlambda} can be represented by the product of $B_\Delta^T$ with a Lagrange multiplier vector $\lambda$. If $\lambda$ is not in $\Lambda$, i.e., not in the range of $B_\Delta$, it can always be written as the sum of its components in $\Lambda$ and in the null of $B_\Delta^T$. Then we just take its component in $\Lambda$ as $\lambda$, which does not change the product $B_\Delta^T \lambda$. %We also note that such a $\lambda \in \Lambda$ is unique, since the difference of any two such $\lambda$s must be both in the null of $B_\Delta^T$ and in the range of $B_\Delta$, which can only be zero.
By multiplying $B_{\Delta,D}$ to both sides of \EQ{aboutlambda} and using Lemma \LA{Bdelta}, we have \EQ{lambda}. ~~
\endproof

\begin{mylemma}
\label{lemma:nullspace}
The basis vector in the null space of \EQ{bigeq},  corresponding to the one-dimensional null space of the original incompressible Stokes system~\EQ{matrix}, is
\be
\label{equation:basis}
\left( \begin{array}{cccccc}
{\bf 0}, & 1_{p_I}, & {\bf 0}, & {\bf 0}, & 1_{p_\Gamma}, & -B_{\Delta,D}[B_{I\Delta}^T ~~
B_{\Gamma\Delta}^T]\left[\begin{array}{c}1_{p_I}\\
    1_{p_\Gamma}\end{array}\right] \end{array}\right).
\ee
\end{mylemma}

\beginproof Since the null space of \EQ{matrix} consists of all constant pressures, substituting the vector \EQ{basis} into \EQ{bigeq} gives zero blocks on the  right-hand side, except at the third block where
\be
\label{equation:rightside}
{\bf f}_{\Delta} = [B_{I\Delta}^T ~~ B_{\Gamma\Delta}^T]\left[\begin{array}{c}1_{p_I}\\
1_{p_\Gamma}\end{array}\right] - B^T_\Delta B_{\Delta,D}[B_{I\Delta}^T ~~ B_{\Gamma\Delta}^T]\left[\begin{array}{c}1_{p_I}\\
   1_{p_\Gamma}\end{array}\right],
\ee
which also equals zero from \EQ{aboutlambda} and \EQ{lambda}  in Lemma \LA{faceIntegral}. \qquad
\endproof

\section{A reduced symmetric positive semi-definite system}
\label{section:Gmatrix}

The system \EQ{bigeq} can be reduced to a Schur complement problem for the variables $\left(p_{\Gamma}, ~\lambda \right)$.
Since the leading four-by-four block of the coefficient matrix in \EQ{bigeq} is invertible, the variables  $\left( {\bf u}_I, ~p_I, ~{\bf u}_{\Delta},
~{\bf u}_{\Pi} \right)$ can be eliminated and we obtain
\begin{equation}
\label{equation:spd}
G \left[ \begin{array}{c}
p_\Gamma         \\[0.8ex]
\lambda
\end{array} \right] ~ = ~ g,
\end{equation}
where
\begin{equation}
\label{equation:Gmatrixgvec}
G = B_C \widetilde{A}^{-1} B_C^T, \qquad \qquad
g = B_C\widetilde{A}^{-1}\left[
\begin{array}{l}
{\bf f}_I        \\[0.8ex]
0                \\[0.8ex]
{\bf f}_{\Delta} \\[0.8ex]
{\bf f}_\Pi
\end{array}
\right],
\end{equation}
with
\begin{equation}
\label{equation:AtildeBc}
~~ \widetilde{A} = \left[
\begin{array}{cccc}
A_{II}      & B_{II}^T       & A_{I \Delta}       & A_{I \Pi}       \\[0.8ex]
B_{II}      & 0              & B_{I \Delta}       & B_{I \Pi}       \\[0.8ex]
A_{\Delta I}& B_{I \Delta} ^T& A_{\Delta\Delta}   & A_{\Delta \Pi}  \\[0.8ex]
A_{\Pi I}   & B_{I \Pi}^T    & A_{\Pi \Delta}     & A_{\Pi \Pi}   \end{array}
\right]
\quad \mbox{and} \quad
B_C=\left[
\begin{array}{cccc}
B_{\Gamma I} & 0 & B_{\Gamma \Delta} & B_{\Gamma \Pi} \\[0.8ex]
0            & 0 & B_{\Delta}        & 0              \end{array}
\right].
\end{equation}

We can see that $-G$ is the Schur complement of the coefficient matrix
of \EQ{bigeq} with respect to the last two row blocks, i.e.,
\[
\left[ \begin{array}{cc} I & 0 \\[0.8ex] -B_C \widetilde{A}^{-1} & I \end{array} \right]
\left[ \begin{array}{cc} \widetilde{A} & B_C^T \\[0.8ex] B_C & 0 \end{array} \right]
\left[ \begin{array}{cc} I & - \widetilde{A}^{-1} B_C^T \\[0.8ex] 0  & I \end{array} \right] =
\left[ \begin{array}{cc} \widetilde{A} & 0 \\[0.8ex] 0 & -G \end{array} \right].
\]
From the Sylvester law of inertia, namely, the number of positive, negative, and zero eigenvalues of a symmetric matrix is invariant under a change of coordinates, we can see that the number of zero eigenvalues of $G$ is the same as the number of zero eigenvalues (with multiplicity counted) of the original coefficient matrix of \EQ{bigeq}, and all other eigenvalues of $G$ are positive. Therefore $G$ is symmetric positive semi-definite. The basis vectors of the null space of $G$ also inherit those from the null space of \EQ{bigeq}, and the only interesting basis vector is
\be
\label{equation:onenull}
\left( \begin{array}{cc} 1_{p_\Gamma}, & - B_{\Delta,D}[B_{I\Delta}^T ~~ B_{\Gamma\Delta}^T]\left[\begin{array}{c}1_{p_I}\\ 1_{p_\Gamma} \end{array} \right]
\end{array} \right),
\ee
which is derived from Lemma \LA{nullspace}. The other null space vectors of $G$ are all vectors with $\lambda$ in the null of $B_\Delta^T$ and $p_\Gamma = 0$. The range of G contain all vectors orthogonal to those null vectors.
Denote $X = Q_\Gamma \bigoplus \Lambda$, where, as defined earlier, $\Lambda$ is the range of $B_\Delta$. Then the range of $G$, denoted by $R_G$, is the subspace of $X$ orthogonal to \EQ{onenull}, i.e.,
\be
\label{equation:Grange}
R_G=\left\{  \left[ \begin{array}{c}
g_{p_\Gamma}         \\[0.8ex]
g_{\lambda}
\end{array} \right] \in X ~ {\Big|} ~ g_{p_\Gamma}^T 1_{{p_\Gamma}} -
g_{\lambda}^T \left(B_{\Delta,D}[B_{I\Delta}^T ~~ B_{\Gamma\Delta}^T]\left[\begin{array}{c}1_{{p_I}}\\ 1_{{p_\Gamma}}\end{array}\right]\right)=0\right\}.
\ee

The restriction of $G$ to its range $R_G$ is positive definite.
The fact that the solution of \EQ{bigeq} always exists for any given
$\left( {\bf f}_I, ~{\bf f}_{\Delta}, ~{\bf f}_{\Pi} \right)$ on the
right-hand side implies that the solution of~\EQ{spd} exits for any
$g$ defined by \EQ{Gmatrixgvec}. Therefore $g \in R_G$. When the conjugate
gradient method (CG) is applied to solve \EQ{spd} with zero initial
guess, all the iterates are in the Krylov subspace generated by $G$
and $g$, which is also a subspace of $R_G$, and where the CG cannot
break down.
After obtaining $\left( p_{\Gamma}, ~\lambda \right)$ from solving \EQ{spd}, the other components $\left( {\bf u}_I, ~p_I, ~{\bf u}_{\Delta}, ~{\bf u}_{\Pi} \right)$ in \EQ{bigeq} are obtained by back substitution.

In the rest of this section, we discuss the implementation of multiplying $G$ by a vector. The main operation is the product of $\Atilde^{-1}$ with a vector, cf. \EQ{Gmatrixgvec}.  We denote
\[ A_{rr} = \left[ \begin{array}{ccc}
A_{II}       & B_{II}^T        & A_{I \Delta}     \\[0.8ex]
B_{II}       & 0               & B_{I \Delta}     \\[0.8ex]
A_{\Delta I} & B_{I \Delta} ^T & A_{\Delta\Delta} \end{array} \right] ,  \quad
A_{\Pi r} = A_{r \Pi}^T = \left[ A_{\Pi I}  \quad B_{I \Pi}^T  \quad A_{\Pi \Delta} \right], \quad  f_r = \left[ \begin{array}{l}
{\bf f}_I        \\[0.8ex]
0                \\[0.8ex]
{\bf f}_{\Delta} \end{array} \right],
\]
and define the Schur complement
\[
S_{\Pi} = A_{\Pi \Pi} - A_{\Pi r} A_{rr}^{-1} A_{r \Pi},
\]
which is symmetric positive definite from the Sylvester law of inertia. $S_\Pi$ defines the coarse level problem in the algorithm.  The product
\[
   \left[
   \begin{array}{cccc}
   A_{II}       & B_{II}^T        & A_{I \Delta}      & A_{I \Pi}      \\[0.8ex]
   B_{II}       & 0               & B_{I \Delta}      & B_{I \Pi}      \\[0.8ex]
   A_{\Delta I} & B_{I \Delta} ^T & A_{\Delta\Delta}  & A_{\Delta \Pi} \\[0.8ex]
   A_{\Pi I}    & B_{I \Pi}^T     & A_{\Pi \Delta}    & A_{\Pi \Pi}
   \end{array}
   \right]^{-1} \left[
   \begin{array}{l}
   {\bf f}_I        \\[0.8ex]
   0                \\[0.8ex]
   {\bf f}_{\Delta} \\[0.8ex]
   {\bf f}_\Pi
   \end{array}
   \right]
\]
can then be represented by
\[
\left[ \begin{array}{c} A_{rr}^{-1} f_r \\[0.8ex] \vvec{0} \end{array} \right] ~ + ~
\left[ \begin{array}{c} -A_{rr}^{-1} A_{r \Pi} \\[0.8ex] I_\Pi \end{array}  \right] ~ S_{\Pi}^{-1} ~
\left({\bf f}_\Pi - A_{\Pi r} A_{rr}^{-1} f_r \right),
\]
which requires solving the coarse level problem once and independent subdomain Stokes problems with Neumann type boundary conditions twice.
%We note that here the coarse level problem matrix $S_\Pi$ is symmetric positive definite, since it is the Schur complement on the coarse level velocity degrees of freedom.

\section{Preliminary results}
\label{section:techniques}
%In this section, we will give some preliminary results needed in our analysis. Many of them have been discussed in \cite{Tu:2012:Stokes} for two dimensions and are still valid for three dimensions.
%
%We first define certain norms for several vector/function spaces.
Denote
\be
\label{equation:Wtilde}
\vvec{\Wtilde} = {\bf W}_I \bigoplus \vvec{\Wtilde}_{\Gamma} = {\bf W}_I \bigoplus {\bf W}_{\Delta} \bigoplus {\bf W}_\Pi.
\ee
For any $\vvec{w}$ in $\vvec{\Wtilde}$, denote its restriction to subdomain $\Omega_i$ by ${\bf w}^{(i)}$. A subdomain-wise $H^1$-seminorm can be defined for functions in $\vvec{\Wtilde}$ by
\[
|\vvec{w}|^2_{H^1} = \sum_{i=1}^N |\vvec{w}^{(i)}|^2_{H^1(\Omega_i)}.
\]

We also define
\[
\Vtilde = {\bf W}_I \bigoplus Q_I \bigoplus {\bf W}_{\Delta} \bigoplus {\bf W}_\Pi,
\]
and its subspace
\be
\label{equation:W0}
\Vtilde_0 = \left\{ v = \left( {\bf w}_I, ~p_I, ~{\bf w}_{\Delta}, ~{\bf w}_{\Pi} \right) \in \Vtilde ~ \big| ~  B_{I I} \vvec{w}_I + B_{I\Delta} \vvec{w}_\Delta + B_{I\Pi} \vvec{w}_\Pi = 0 \right\}.
\ee
For any $v = \left( {\bf w}_I, ~p_I, ~{\bf w}_{\Delta}, ~{\bf w}_{\Pi} \right) \in \Vtilde_0$, let
$\vvec{w} =  \left( {\bf w}_I, ~{\bf w}_{\Delta}, ~{\bf w}_{\Pi} \right) \in \vvec{\Wtilde}$. Then
\begin{eqnarray}\label{equation:wg0}
\left< v, v \right>_{\widetilde{A}} & = &
\left[ \begin{array}{l} {\bf w}_I \\[0.8ex] {\bf w}_{\Delta} \\[0.8ex] {\bf w}_\Pi \end{array} \right]^T
\left[ \begin{array}{ccc}
A_{II}        & A_{I \Delta}      & A_{I \Pi}      \\[0.8ex]
A_{\Delta I}  & A_{\Delta\Delta}  & A_{\Delta \Pi} \\[0.8ex]
A_{\Pi I}     & A_{\Pi \Delta}    & A_{\Pi \Pi} \end{array} \right]
\left[ \begin{array}{l} {\bf w}_I \\[0.8ex] {\bf w}_{\Delta} \\[0.8ex] {\bf w}_\Pi \end{array} \right] \nonumber \\[0.8ex]
& = & \sum_{i=1}^N \left[ \begin{array}{c} {\bf w}_I^{(i)} \\ {\bf w}_{\Delta}^{(i)} \\ {\bf w}_{\Pi}^{(i)} \end{array} \right]^T
\left[ \begin{array}{cccc}
A_{II}^{(i)}       & A_{I \Delta}^{(i)}      & A_{I \Pi}^{(i)}      \\[0.8ex]
A_{\Delta I}^{(i)} & A_{\Delta\Delta}^{(i)}  & A_{\Delta \Pi}^{(i)} \\[0.8ex]
A_{\Pi I}^{(i)}    & A_{\Pi \Delta}^{(i)}    & A_{\Pi \Pi}^{(i)}
\end{array} \right] \left[ \begin{array}{c} {\bf w}_I^{(i)} \\ {\bf w}_{\Delta}^{(i)} \\ {\bf w}_{\Pi}^{(i)} \end{array} \right]
= \sum_{i=1}^N \left| \left[ \begin{array}{c} {\bf w}_I^{(i)} \\ {\bf w}_{\Delta}^{(i)} \\ {\bf w}_{\Pi}^{(i)} \end{array} \right] \right|_{H^1(\Omega_i)}^2 \label{equation:W0n} \\[0.8ex]
& = & |\vvec{w}|^2_{H^1},  \nonumber
\end{eqnarray}
where the superscript ${}^{(i)}$ is used to represent the restrictions
of corresponding vectors and matrices to subdomain $\Omega_i$. We can
see from \EQ{W0n} that for any $v\in \Vtilde_0$, the value $\left< v,
  v\right>_{\widetilde{A}}$ is independent of its pressure component
$p_I$. $\left< \cdot, \cdot \right>_{\widetilde{A}}$ defines a semi-inner
product on $\Vtilde_0$;  $\left< v, v \right>_{\widetilde{A}}=0$  if
and only if the velocity component of $v$ is constant on $\Omega$ and is in fact zero due to the zero boundary condition on $\partial \Omega$, while its pressure component can be arbitrary.

Denote
\begin{equation}\label{equation:Btilde}
\widetilde{B} = \left[ \begin{array}{ccc}
B_{II}      & B_{I \Delta}       & B_{I \Pi}      \\[0.8ex]
B_{\Gamma I}& B_{\Gamma \Delta}  & B_{\Gamma \Pi}
\end{array} \right],
\end{equation}
cf. \EQ{bigeq}. The following lemma on the stability of $\widetilde{B}$ can be found at \cite[Lemma 5.1]{li12}.

\begin{mylemma}
\label{lemma:BtildeStability}
For any $\vvec{w} \in \vvec{\Wtilde}$ and $q \in Q$, $\left< \Btilde {\bf w}, q \right> \leq | \vvec{w} |_{H^1} \| q \|_{L^2}$.
\end{mylemma}

%\beginproof \begin{eqnarray*}
%\left<\Btilde {\bf w}, q \right>^2 & = & \left( \sum_{i=1}^N \int_{\Omega_i} \nabla\cdot {\bf w}^{(i)} q \right)^2\leq \left( \sum_{i=1}^N \sqrt{\int_{\Omega_i} | \nabla {\bf w}^{(i)} |^2} \sqrt{\int_{\Omega_i} q^2} \right)^2  \\[0.8ex]
%& \le & \left( \sum_{i=1}^N \int_{\Omega_i} | \nabla {\bf w}^{(i)} |^2 \right) \left( \sum_{i=1}^N \int_{\Omega_i} q^2 \right) = | \vvec{w} |^2_{H^1} \| q \|^2_{L^2}. \qquad \Box
%\end{eqnarray*}

The following lemma will also be used and can be found at \cite[Lemma~2.3]{paulo2003}.
\begin{mylemma}
\label{lemma:paul}
Let  $(\vvec{u}, p) \in \vvec{W} \bigoplus Q$ satisfy
\be \left[
\begin{array}{cc}
A       & B^T      \\[0.8ex]
B       & 0
\end{array}
\right] \left[
\begin{array}{l}
\vvec{u}        \\[0.8ex]
p
\end{array}
\right]
=\left[
\begin{array}{l}
\vvec{f}        \\[0.8ex]
g
\end{array}
\right],
\ee
where $A$ and $B$ are as in \EQ{matrix}, $\vvec{f} \in \vvec{W}$, and $g \in Im(B) \subset Q$. Let $\beta$ be the inf-sup constant specified in \EQ{infsupMatrix}.
Then
\[
\| \vvec{u} \|_A \le \| \vvec{f} \|_{A^{-1}} + \frac{1}{\beta} \| g \|_{Z^{-1}},
\]
where $Z$ is the mass matrix defined in Section \ref{section:FEM}.
\end{mylemma}

\section{Jump operators and preconditioners}
\label{section:jump}

We first define certain jump operators across the subdomain interface $\Gamma$, which will be used for the analysis of the preconditioners.

Denote the restriction operator from $\Vtilde$ onto ${\bf W}_{\Delta}$ by $\widetilde{R}_\Delta$, i.e., for any $v = \left( {\bf w}_I, ~p_I, ~{\bf w}_{\Delta}, ~{\bf w}_{\Pi} \right) \in \Vtilde$,
$\widetilde{R}_\Delta v = {\bf w}_{\Delta}$. Define $P_{D,L}: \Vtilde \rightarrow \Vtilde$, by
\[
P_{D,L} = \widetilde{R}_{\Delta}^T B_{\Delta, D}^T B_{\Delta} \widetilde{R}_\Delta.
\]
Following this definition, given any $v = \left( {\bf w}_I, ~p_I, ~{\bf w}_{\Delta}, ~{\bf w}_{\Pi} \right) \in \Vtilde$, the dual velocity component of $P_{D,L} v$, on any subdomain interface node $x$ in subdomain $\Omega_i$, is given by, cf.~\cite[Equation (6.70)]{Toselli:2004:DDM},
\[
\left( \widetilde{R}_\Delta \left( P_{D,L} v \right) \right)^{(i)} (x) = \sum_{j \in \N_{x}} \delta^{\dagger}(x) \left( {\bf w}_\Delta^{(i)} (x) - {\bf w}_\Delta^{(j)} (x) \right),
\]
which represents the so-called jump of the dual velocity component ${\bf w}_{\Delta}$ across the subdomain interface $\Gamma$. All other components of $P_{D,L} v$ equal zero. We also have
\begin{eqnarray}\label{equation:PD1}
\left< P_{D,L} v, P_{D,L} v \right>_{\widetilde{A}} &=&
(\widetilde{R}_{\Delta}^T B_{\Delta, D}^T B_{\Delta}
\widetilde{R}_\Delta v)^T\widetilde{A}(\widetilde{R}_{\Delta}^T B_{\Delta, D}^T
B_{\Delta} \widetilde{R}_\Delta v)\nonumber\\
&=&\left< B_{\Delta, D}^T B_{\Delta} {\bf w}_{\Delta}, B_{\Delta, D}^T B_{\Delta} {\bf w}_{\Delta} \right>_{A_{\Delta \Delta}}.
\end{eqnarray}
Together with \EQ{W0n}, we have the following lemma, which can be found at \cite[Section~6.1]{li07}.

\begin{mylemma}
\label{lemma:jump}
There exists a constant $C$ and a function $\Phi_L(H/h)$, such that for all $v \in \Vtilde_0$,
$\left< P_{D,L}  v, P_{D,L}  v \right>_{\widetilde{A}} \leq  C \Phi_L(H/h) \left< v, v \right>_{\widetilde{A}}$. Here, $\Phi_L(H/h) = (H/h) (1 +\log{(H/h)})$,  when the coarse level space is spanned by the subdomain vertex nodal basis functions and subdomain edge-cutoff functions corresponding to each velocity component.
\end{mylemma}

When applying $P_{D,L}$ to a vector, the jump of the dual subdomain
interface velocities is extended by zero to the interior of
subdomains. To improve the stability of the jump operator, the jump can be extended to the interior of subdomains by subdomain discrete harmonic extension. We define a Schur complement operator $H^{(i)}_{\Delta}: {\bf W}^{(i)}_{\Delta} \rightarrow {\bf W}^{(i)}_{\Delta}$ by, for any ${\bf u}_{\Delta}^{(i)} \in {\bf W}^{(i)}_{\Delta}$,
\begin{equation}
\label{equation:SiDelta}
\left[ \begin{array}{cc}
A_{II}^{(i)}            & A_{I \Delta}^{(i)}     \\[0.8ex]
%B_{II}^{(i)}      & 0                   & B_{I \Delta}^{(i)}     \\[0.8ex]
A_{\Delta I}^{(i)} & A_{\Delta\Delta}^{(i)}
\end{array} \right] \left[ \begin{array}{c}
{\bf u}_I^{(i)}        \\[0.8ex]
{\bf u}_{\Delta}^{(i)}
\end{array} \right] =
\left[ \begin{array}{l}
{\bf 0}          \\[0.8ex]
%0                \\[0.8ex]
H_{\Delta}^{(i)} {\bf u}_\Delta^{(i)}
\end{array} \right] \mbox{ . }
\end{equation}
To multiply $H_{\Delta}^{(i)}$ by a vector $\vvec{u}_\Delta^{(i)}$, a subdomain elliptic problem on $\Omega_i$ with given boundary velocity $\vvec{u}_\Delta^{(i)}$ and $\vvec{u}^{(i)}_\Pi=\vvec{0}$ needs to be solved. We let $H_{\Delta}: {\bf W}_{\Delta} \rightarrow {\bf W}_{\Delta}$ to represent the direct sum of  $H^{(i)}_\Delta, i = 1, \ldots, N$.

Using $H^{(i)}_{\Delta}$, we define the second jump operator $P_{D,D}: \Vtilde \rightarrow \Vtilde$, by: for any given $v = \left( {\bf w}_I, ~p_I, ~{\bf w}_{\Delta}, ~{\bf w}_{\Pi} \right) \in \Vtilde$, the subdomain interior velocity part of $P_{D,D} v$ on each subdomain $\Omega_i$ is taken as ${\bf u}_I^{(i)}$ in the solution of \EQ{SiDelta}, with given subdomain boundary velocity ${\bf u}_{\Delta}^{(i)} = B_{\Delta, D}^{(i)^T} B_{\Delta} {\bf w}_{\Delta}$. Here $B_{\Delta, D}^{(i)^T}$ represents restriction of $B_{\Delta, D}^T$ on subdomain $\Omega_i$ and is a map from $\Lambda$ to $\vvec{W}^{(i)}_\Delta$. The other components of $P_{D,D} v$ are kept zero. Therefore
\begin{eqnarray}
\qquad & & \left<P_{D,D} v,P_{D,D} v\right>_{\Atilde} = \sum_{i=1}^N
\left[ \begin{array}{c}
{\bf u}_I^{(i)}        \\[0.8ex]
%p_I^{(i)}              \\[0.8ex]
{\bf u}_{\Delta}^{(i)}
\end{array} \right]^T  \left[ \begin{array}{ccc}
A_{II}^{(i)}          & A_{I \Delta}^{(i)}     \\[0.8ex]
%B_{II}^{(i)}                      & B_{I \Delta}^{(i)}     \\[0.8ex]
A_{\Delta I}^{(i)}& A_{\Delta\Delta}^{(i)}
\end{array} \right] \left[ \begin{array}{c}
{\bf u}_I^{(i)}        \\[0.8ex]
{\bf u}_{\Delta}^{(i)}
\end{array} \right]  \nonumber \\
& = &  \sum_{i=1}^N  {\bf u}_\Delta^{(i)^T} H_{\Delta}^{(i)} {\bf u}_\Delta^{(i)} = \sum_{i=1}^N {\bf w}_{\Delta}^T B_{\Delta}^T
B_{\Delta, D}^{(i)} H_{\Delta}^{(i)} B_{\Delta, D}^{(i)^T} B_{\Delta}
{\bf w}_{\Delta} \nonumber\\
& = & \sum_{i=1}^N \left|\left[\begin{array}{c} B_{\Delta,
        D}^{(i)^T} B_{\Delta} {\bf w}_{\Delta}\\
      0 \end{array}\right]\right|_{H^{1/2}(\partial
  \Omega^i)}^2 \leq C
\Phi_D(H/h) \sum_{i=1}^N \left|\left[ \begin{array}{c} {\bf w}_{\Delta}^{(i)} \\ {\bf w}_{\Pi}^{(i)} \end{array} \right] \right|^2_{H^{1/2}(\partial
  \Omega^i)}\label{equation:PD2}\\
&\leq&
C \Phi_D(H/h) \sum_{i=1}^N \left| \left[ \begin{array}{c} {\bf
        w}_I^{(i)} \\ {\bf w}_{\Delta}^{(i)} \\ {\bf
        w}_{\Pi}^{(i)} \end{array} \right]
\right|_{H^1(\Omega_i)}^2= C \Phi_D(H/h) |{\bf w}|_{H^1(\Omega_i)}^2. \nonumber
\end{eqnarray}
The first inequality in \EQ{PD2} is a well established
result, cf., \cite[Lemma 6.36]{Toselli:2004:DDM}. Since for
any $v\in \Vtilde_0$, $\left< v, v
\right>_{\widetilde{A}}=|{\bf w}|_{H^1(\Omega_i)}^2$, cf.~\EQ{W0n},  we have the following lemma.

\begin{mylemma}
\label{lemma:jumptwo} There exists a constant $C$ and a function $\Phi_D(H/h)$, such that for all $v \in \Vtilde_0$, $\left< P_{D,D} v, P_{D,D} v \right>_{\widetilde{A}}\leq C \Phi_D(H/h) \left< v, v \right>_{\widetilde{A}}$.  Here $\Phi_D(H/h) = (1 + \log{(H/h)})^2$,  when the coarse level space is spanned by the subdomain vertex nodal basis functions and subdomain edge-cutoff functions corresponding to each velocity component.
\end{mylemma}

To introduce the preconditioners, we write $G$, defined in \EQ{Gmatrixgvec} and \EQ{AtildeBc}, in a two-by-two block structure.  Denote the first row of $B_C$ by
\[
\widetilde{B}_{\Gamma} =
\left[ B_{\Gamma I} \quad 0  \quad B_{\Gamma \Delta}  \quad  B_{\Gamma \Pi} \right],
\]
and note that $\widetilde{R}_{\Delta}$ is the restriction operator from $\Vtilde$ onto ${\bf W}_{\Delta}$.  Then $G$ can be written as
\begin{equation}
\label{equation:Gtwo}
G = \left[ \begin{array}{cc} G_{p_\Gamma p_\Gamma} & G_{p_\Gamma \lambda} \\[0.8ex] G_{\lambda p_\Gamma} & G_{\lambda \lambda} \end{array} \right],
\end{equation}
where
\begin{eqnarray*}
& G_{p_\Gamma p_\Gamma} = \widetilde{B}_{\Gamma} \widetilde{A}^{-1} \widetilde{B}_{\Gamma}^T, \qquad
G_{p_\Gamma \lambda} = \widetilde{B}_{\Gamma} \widetilde{A}^{-1} \widetilde{R}_{\Delta}^T B_{\Delta}^T, & \\ [0.8ex]
& G_{\lambda p_\Gamma} = B_{\Delta} \widetilde{R}_\Delta \widetilde{A}^{-1} \widetilde{B}_{\Gamma}^T, \qquad
G_{\lambda \lambda} = B_{\Delta} \widetilde{R}_\Delta \widetilde{A}^{-1} \widetilde{R}_{\Delta}^T B_{\Delta}^T. &
\end{eqnarray*}

We consider a block diagonal preconditioner for~\EQ{spd}. As for two-dimensional problems, the first diagonal block $G_{p_\Gamma p_\Gamma}$ of $G$ can be shown spectrally equivalent to
$h^3 I_{p_\Gamma}$, where $I_{p_\Gamma}$ is the identity matrix of the
same dimension as $G_{p_\Gamma p_\Gamma}$; see
\cite{li12,Tu:2012:Stokes}. Therefore, in the following block diagonal
preconditioners, the inverse of $G_{p_\Gamma p_\Gamma}$ is
approximated by $\alpha h^{-3} I_{p_\Gamma}$. Here $\myalpha$ is a given constant. We will
show in the next section that $\alpha$ has only a minor effect on the condition number bound of the preconditioned operator and its value is typically taken as 1, cf. Remark \ref{remark:alpha}. We introduce $\alpha$ in the preconditioner just for the convenience in the numerical experiments to demonstrate the convergence rates of the proposed algorithm.

The inverse of the second diagonal block $B_{\Delta} \widetilde{R}_\Delta \widetilde{A}^{-1} \widetilde{R}_{\Delta}^T B_{\Delta}^T$, can be approximated by the lumped block
\begin{equation}\label{equation:precondlumped}
M^{-1}_{\lambda,L} = B_{\Delta, D} \widetilde{R}_\Delta \widetilde{A} \widetilde{R}_{\Delta}^T B_{\Delta, D}^T.
\end{equation}
This leads to the following lumped preconditioner for solving \EQ{spd}
\be
\label{equation:lumped}
M_L^{-1} = \left[ \begin{array}{cc} \alpha h^{-3} I_{p_\Gamma} & \\[0.8ex]
& M^{-1}_{\lambda,L}  \end{array} \right].
\ee
Applying subdomain discrete harmonic extensions in the preconditioning step, we have the following Dirichlet preconditioner
\be
\label{equation:Dirichlet}
M_{D}^{-1} = \left[ \begin{array}{cc} \alpha h^{-3} I_{p_\Gamma} & \\[0.8ex]
& M^{-1}_{\lambda,D}  \end{array} \right],
\ee
where
\be
\label{equation:subdomainwise}
M^{-1}_{\lambda,D} =B_{\Delta, D} H_{\Delta} B_{\Delta, D}^T.
\ee

We can see from Lemma \LA{BdeltaNull} that both $M^{-1}_{\lambda,L}$ and $M^{-1}_{\lambda,D}$ are symmetric positive definite when restricted on $\Lambda$. Therefore both the lumped and the Dirichlet preconditioners $M^{-1}_L$ and $M^{-1}_D$ are symmetric positive definite in the range of $G$.

\section{Condition number bounds}
\label{section:convergence}

In the following, we use the same framework to establish the condition number bounds for both
lumped and Dirichlet preconditioned operators $M_L^{-1} G$ and
$M_D^{-1} G$. Let $M^{-1}$, $M^{-1}_\lambda$, $P_D$, and $\Phi$ to
represent both $M^{-1}_L$, $M^{-1}_{\lambda,L}$, $P_{D, L}$, $\Phi_L$, for the lumped preconditioner case, and $M^{-1}_D$, $M^{-1}_{\lambda,D}$, $P_{D, D}$, $\Phi_D$, for the Dirichlet preconditioner case, respectively, when they apply in the proofs.

When the conjugate gradient method is applied to solving the preconditioned system
\begin{equation}
\label{equation:Mspd}
M^{-1} G x ~ = ~ M^{-1} g,
\end{equation}
with zero initial guess, all iterates belong to the Krylov subspace generated by the operator $M^{-1} G$ and the vector $M^{-1} g$, which is a subspace of the range of $M^{-1} G$. We denote the range of $M^{-1} G$ by $R_{M^{-1} G}$ and note that both preconditioners are symmetric positive definite in the range of $G$. We have the following lemma, cf.~\cite[Lemma 6]{Tu:2012:Stokes}.

\begin{mylemma}
\label{lemma:CG}
The conjugate gradient method applied to solving \EQ{Mspd} with zero initial guess cannot break down.
\end{mylemma}

\beginproof We just need to show that for any $0 \neq x \in R_{M^{-1} G}$, $\left<x, G x \right> \neq 0$, i.e., to show $G x \neq 0$. Let $0 \neq x = M^{-1} G y$, for a certain $y \in X$ and $y \neq 0$. Then $G x = G M^{-1} G y$, which cannot be zero since $G y \neq 0$ and $y^T  G M^{-1} G y \neq 0$. \qquad
\endproof

The following lemma will be used to provide the upper eigenvalue bound of the preconditioned operator. It is similar to \cite[Lemma 6.4]{li12} and \cite[Lemmas 8 and 11]{Tu:2012:Stokes}.

\begin{mylemma}\label{lemma:upper}
There exists a constant $C$, such that for all $v \in \Vtilde_0$,
\[
\left<M^{-1} B_C v, B_C v \right>\le C\left(\myalpha+\Phi(H/h)\right)\left<\Atilde v, v \right>,
\]
where $\Phi(H/h)$ is defined in Lemmas \LA{jump} and \LA{jumptwo}, respectively.
\end{mylemma}

\beginproof Given $v = \left( {\bf w}_I, ~q_I, ~{\bf w}_{\Delta}, ~{\bf w}_{\Pi}
\right) \in \Vtilde_0$, let $g_{p_\Gamma} = B_{\Gamma I} {\bf w}_I + B_{\Gamma \Delta} {\bf w}_\Delta + B_{\Gamma\Pi} {\bf w}_\Pi$. From \EQ{AtildeBc}, \EQ{precondlumped}--\EQ{subdomainwise}, \EQ{PD1}, and \EQ{PD2}, we have
\BA
\label{equation:MBW}
\left<M^{-1}B_Cv,B_Cv\right> & = & \alpha h^{-3} \left<g_{p_\Gamma}, g_{p_\Gamma} \right>+
\left(B_\Delta\Rtilde_\Delta v\right)^T M^{-1}_{\lambda} B_\Delta\Rtilde_\Delta v\nonumber\\
&=& \alpha h^{-3} \left< g_{p_\Gamma}, g_{p_\Gamma} \right>+\left<P_D v,P_D v\right>_\Atilde\nonumber\\
&\le& \alpha h^{-3} \left< g_{p_\Gamma}, g_{p_\Gamma} \right>+ C \Phi(H/h)\left<v,v\right>_\Atilde,
\EA
where we used Lemmas \LA{jump} and \LA{jumptwo} for the last inequality. It is  sufficient to bound the first term of the right-hand side in the above inequality.

Since $v \in \Vtilde_0$, we have $B_{I I} \vvec{w}_I + B_{I\Delta} \vvec{w}_\Delta + B_{I\Pi} \vvec{w}_\Pi = 0$, cf.~\EQ{W0}. Then
\begin{eqnarray*}
\left<g_{p_\Gamma}, g_{p_\Gamma}\right> & = & \left[ \begin{array}{c} B_{I I} \vvec{w}_I + B_{I\Delta} \vvec{w}_\Delta + B_{I\Pi} \vvec{w}_\Pi \\ B_{\Gamma I} {\bf w}_I + B_{\Gamma \Delta} {\bf w}_\Delta + B_{\Gamma\Pi} {\bf w}_\Pi \end{array} \right]^T
\left[ \begin{array}{c} B_{I I} \vvec{w}_I + B_{I\Delta} \vvec{w}_\Delta + B_{I\Pi} \vvec{w}_\Pi \\ B_{\Gamma I} {\bf w}_I + B_{\Gamma \Delta} {\bf w}_\Delta + B_{\Gamma\Pi} {\bf w}_\Pi \end{array} \right] \\
& = & \left<\Btilde {\bf w},  \Btilde {\bf w} \right>,
\end{eqnarray*}
where $\Btilde$ is defined in \EQ{Btilde} and $\vvec{w} = \left( {\bf w}_I, ~{\bf w}_{\Delta}, ~{\bf w}_{\Pi} \right) \in \vvec{\Wtilde}$. From \EQ{massmatrix} and the stability of $\Btilde$, cf. Lemma \LA{BtildeStability}, we have
\begin{eqnarray}
~~~ h^{-3} \left<g_{p_\Gamma}, g_{p_\Gamma}\right> & = & h^{-3} \left<\Btilde {\bf w},  \Btilde {\bf w} \right> \leq C \left<\Btilde {\bf w},  \Btilde {\bf w} \right>_{Z^{-1}} =  C \max_{q \in Q} \frac{\left<\Btilde {\bf w}, q \right>^2}{\left<q, q\right>_Z} \label{equation:boundBu}\\
& \le & C \max_{q \in Q} \frac{| \vvec{w} |^2_{H^1}  \| q \|^2_{L^2}}{\| q \|^2_{L^2}} = C | \vvec{w} |^2_{H^1} = C \left<v, v\right>_\Atilde, \nonumber
\end{eqnarray}
where for the last equality, we used \EQ{W0n}. \qquad
\endproof

The following lemma will be used to provide the lower eigenvalue bound of the preconditioned operator. In \cite[Lemmas 9 and 12]{Tu:2012:Stokes}, the lower eigenvalue bounds for the lumped and Dirichlet preconditioners were analyzed differently. In the analysis of the Dirichlet preconditioner, subdomain discrete Stokes extensions were used. Such extensions require enforcing the same type divergence free subdomain boundary velocity conditions as discussed in \cite{li06}, even though they are not necessary for implementing the algorithm in~\cite{Tu:2012:Stokes}. The new proof given in the next lemma works for both lumped and Dirichlet preconditioners. It does not use the subdomain Stokes extensions and those additional subdomain divergence free boundary conditions are no longer needed. For both type of preconditioners, the coarse level velocity space can be chosen as simple as for solving scalar elliptic problems corresponding to each velocity component.

\begin{mylemma}
\label{lemma:lower}
There exists a constant $C$, such that for any nonzero $y = (g_{p_{\Gamma}}, g_\lambda) \in R_G$, there exits $v \in \Vtilde_0$, which satisfies $B_C v = y$, $\left <v, v\right>_\Atilde \neq 0$, and

\centerline{$\left <\Atilde v, v\right> \le C\max\left\{ 1,\frac{1}{\myalpha} \right\} \left(1+\frac{1}{\beta^2}\right)  \left< M^{-1}y, y \right>$.}
\end{mylemma}

\beginproof Given $y = (g_{p_{\Gamma}}, g_\lambda) \in R_G$, take ${\bf u}_{\Delta}^{(I)} = B_{\Delta, D}^T g_\lambda$, ${\bf u}_\Pi^{(I)}={\bf 0}$, and $p^{(I)}=0$.
On each subdomain $\Omega_i$, let ${\bf u}_{I}^{(I, i)}$ be zero for the
lumped preconditioner, and be obtained for the Dirichlet preconditioner
through the solution of~\EQ{SiDelta} with given subdomain boundary
values ${\u}_{\Delta}^{(i)} ={\u}_{\Delta}^{(I, i)}$. Let $v^{(I,i)} = \left(\u^{(I,i)}_I, ~p_I^{(I,i)}, ~{\bf u}_{\Delta}^{(I,i)}, ~ {\bf u}^{(I,i)}_\Pi\right)$, the corresponding global vectors
$v^{(I)} = \left(\u^{(I)}_I, ~p_I^{(I)}, ~{\bf u}_{\Delta}^{(I)}, ~ {\bf u}^{(I)}_\Pi\right)$, and $\u^{(I)} = \left(\u^{(I)}_I,  ~{\bf u}_{\Delta}^{(I)}, ~ {\bf u}^{(I)}_\Pi\right)$. Then we have
\be
\label{equation:bcID1}
B_C v^{(I)} = \left[
\begin{array}{cccc}
B_{\Gamma I} & 0 & B_{\Gamma \Delta} & B_{\Gamma \Pi} \\[0.8ex]
0            & 0 & B_{\Delta}        & 0              \end{array}
\right] \left[ \begin{array}{c}
{\bf u}_I^{(I)}        \\[0.8ex]
p_I^{(I)}              \\[0.8ex]
{\bf u}_{\Delta}^{(I)} \\[0.8ex]
{\bf u}_{\Pi}^{(I)}
\end{array} \right] = \left[ \begin{array}{c}
B_{\Gamma I} {\bf u}_I^{(I)} + B_{\Gamma \Delta} {\bf u}_{\Delta}^{(I)} + B_{\Gamma \Pi} {\bf u}_{\Pi}^{(I)}    \\[0.8ex] g_\lambda \end{array} \right],
\ee
where we have used Lemma \LA{Bdelta}.  Also
\begin{eqnarray}
\label{equation:Iprop}
| \u^{(I)} |^2_{H^1} &=& \left[ \begin{array}{c}
{\bf u}_I^{(I)}        \\[0.8ex]
{\bf u}_{\Delta}^{(I)} \\[0.8ex]
{\bf u}_{\Pi}^{(I)}
\end{array} \right]^T \left[
\begin{array}{ccccc}
A_{II}        & A_{I \Delta}      &A_{I \Pi}      \\[0.8ex] A_{\Delta I}  & A_{\Delta\Delta}  &A_{\Delta \Pi} \\[0.8ex]
A_{\Pi I}     & A_{\Pi \Delta}    &A_{\Pi \Pi}
\end{array}
\right] \left[ \begin{array}{c}
{\bf u}_I^{(I)}        \\[0.8ex]
{\bf u}_{\Delta}^{(I)} \\[0.8ex]
{\bf u}_{\Pi}^{(I)}
\end{array} \right] \\\nonumber\\
& =& \left\{ \begin{array}{cc}
|  {\bf u}_{\Delta}^{(I)} |^2_{A_{\Delta\Delta}}, &
\mbox{for lumped preconditioner,}\\[1.0ex]
| {\bf u}_{\Delta}^{(I)} |^2_{H_\Delta}, &
\mbox{for Dirichlet  preconditioner}.\nonumber\end{array}\right.
\end{eqnarray}

We consider a solution to the following fully assembled system of linear equations of the form~\EQ{matrix}: find
$\left({\bf u}_I^{(II)}, ~p_I^{(II)}, ~{\bf u}_{\Gamma}^{(II)}, ~p_\Gamma^{(II)}\right) \in {\bf W}_I \bigoplus Q_I \bigoplus {\bf W}_{\Gamma} \bigoplus Q_\Gamma$, such that
\be
\label{equation:fassembled}
\left[
\begin{array}{cccc}
A_{II}      & B_{II}^T       & A_{I \Gamma}      & B_{\Gamma I}^T     \\[0.8ex]
B_{II}      & 0              & B_{I \Gamma}      & 0                  \\[0.8ex]
A_{\Gamma I}& B_{I \Gamma} ^T& A_{\Gamma\Gamma}  & B_{\Gamma \Gamma}^T \\[0.8ex]
B_{\Gamma I}& 0              & B_{\Gamma \Gamma} & 0
\end{array}
\right]
\left[ \begin{array}{c}
{\bf u}_I^{(II)}        \\[0.8ex]
p_I^{(II)}              \\[0.8ex]
{\bf u}_{\Gamma}^{(II)} \\[0.8ex]
p_{\Gamma}^{(II)}
\end{array} \right] =
\left[ \begin{array}{l}
{\bf 0}        \\[0.8ex]
- B_{I I} {\bf u}_I^{(I)} -  B_{I \Delta} {\bf u}_{\Delta}^{(I)} - B_{I \Pi}  {\bf u}_{\Pi}^{(I)}          \\[0.8ex]
{\bf 0}        \\[0.8ex]
g_{p_{\Gamma}} - B_{\Gamma I} {\bf u}_I^{(I)} -  B_{\Gamma \Delta} {\bf u}_{\Delta}^{(I)} - B_{\Gamma \Pi}  {\bf u}_{\Pi}^{(I)}
\end{array} \right] \mbox{ , }
\ee
where %, from \EQ{schurplus},
 we know that the particularly chosen right-hand side is essentially
\be
\label{equation:rightsidevector}
\left[ \begin{array}{l}
{\bf 0}        \\[0.8ex]
- B_{I I} {\bf u}_I^{(I)} -  B_{I \Delta} {\bf u}_{\Delta}^{(I)}  \\[0.8ex]
{\bf 0}        \\[0.8ex]
g_{p_{\Gamma}} - B_{\Gamma I} {\bf u}_I^{(I)} - B_{\Gamma \Delta} {\bf u}_{\Delta}^{(I)}
\end{array} \right].
\ee
Since $(g_{p_{\Gamma}}, g_\lambda) \in R_G$, we have, cf. \EQ{Grange},
\[
( -B_{I\Delta}{\bf u}^{(I)}_\Delta )^T 1_{p_I}
+ ( g_{p_{\Gamma}}-B_{\Gamma\Delta}{\bf u}^{(I)}_\Delta )^T  1_{p_\Gamma}  =  g_{p_{\Gamma}}^T 1_{p_\Gamma} - g_\lambda^T B_{\Delta, D} \left( B_{I\Delta}^T 1_{p_I} + B_{\Gamma\Delta}^T 1_{p_\Gamma} \right) = 0.
\]
Meanwhile,
\[
(-B_{I I}{\bf u}^{(I)}_I)^T 1_{p_I} + (-B_{\Gamma I}{\bf u}^{(I)}_I)^T  1_{p_\Gamma}
= - \int_\Omega \left( \nabla \cdot {\bf u}^{(I)}_I \right) 1  = 0.
\]
We have that the right-hand side vector \EQ{rightsidevector} has zero average, which implies existence of the solution to \EQ{fassembled}.

Denote ${\bf u}^{(II)} = \left( {\bf u}_I^{(II)}, ~{\bf u}_{\Gamma}^{(II)} \right)$. Then from Lemma \ref{lemma:paul} and \EQ{massmatrix}, we have
\begin{eqnarray}
& & | {\bf u}^{(II)} |^2_{H^1} \leq \frac{1}{\beta^2} \left\| \left[ \begin{array}{l}
- B_{I I} {\bf u}_I^{(I)} -  B_{I \Delta} {\bf u}_{\Delta}^{(I)} - B_{I \Pi}  {\bf u}_{\Pi}^{(I)} \\[0.8ex] g_{p_\Gamma} - B_{\Gamma I} {\bf u}_I^{(I)} -  B_{\Gamma \Delta} {\bf u}_{\Delta}^{(I)} - B_{\Gamma \Pi}  {\bf u}_{\Pi}^{(I)} \end{array} \right] \right\|^2_{Z^{-1}} \nonumber \\
& \le & \frac{1}{\beta^2} \left\| \left[ \begin{array}{l}
B_{I I} {\bf u}_I^{(I)} + B_{I \Delta} {\bf u}_{\Delta}^{(I)} + B_{I \Pi}  {\bf u}_{\Pi}^{(I)} \\[0.8ex] B_{\Gamma I} {\bf u}_I^{(I)} + B_{\Gamma \Delta} {\bf u}_{\Delta}^{(I)} + B_{\Gamma \Pi}  {\bf u}_{\Pi}^{(I)} \end{array} \right] \right\|^2_{Z^{-1}}  +
\frac{1}{\beta^2} \left\| \left[ \begin{array}{l}
0  \\[0.8ex] g_{p_\Gamma}  \end{array} \right] \right\|^2_{Z^{-1}} \nonumber \\
& \leq & \frac{1}{\beta^2} | \u^{(I)} |^2_{H^1}  +  \frac{C}{\beta^2 h^3}\left<g_{p_\Gamma},g_{p_\Gamma}\right>, \label{equation:uIboundD}
\end{eqnarray}
where the bound on the first term is obtained in the same way as in \EQ{boundBu}.

Split the continuous subdomain interface velocity ${\bf u}_{\Gamma}^{(II)}$ into the dual part ${\bf u}_{\Delta}^{(II)}$ and the primal part ${\bf u}_{\Pi}^{(II)}$, and denote $v^{(II)} = \left({\bf u}_I^{(II)}, ~p_I^{(II)}, ~{\bf u}_{\Delta}^{(II)}, ~{\bf u}_{\Pi}^{(II)}\right)$. Let $v = v^{(I)} + v^{(II)}$. Then we have from \EQ{fassembled} that $v \in \widetilde{V}_{0}$, and
\begin{eqnarray*}
B_C v^{(II)} & = & \left[
\begin{array}{cccc}
B_{\Gamma I} & 0 & B_{\Gamma \Delta} & B_{\Gamma \Pi} \\[0.8ex]
0            & 0 & B_{\Delta}        & 0              \end{array}
\right] \left[ \begin{array}{c}
{\bf u}_I^{(II)}        \\[0.8ex]
p_I^{(II)}              \\[0.8ex]
{\bf u}_{\Delta}^{(II)} \\[0.8ex]
{\bf u}_{\Pi}^{(II)}
\end{array} \right] \\
& = & \left[ \begin{array}{c}
g_{p_\Gamma}  - B_{\Gamma I} {\bf u}_I^{(I)} -  B_{\Gamma \Delta} {\bf u}_{\Delta}^{(I)} - B_{\Gamma \Pi}  {\bf u}_{\Pi}^{(I)}    \\[0.8ex] 0  \end{array} \right].
\end{eqnarray*}
Together with \EQ{bcID1}, we have $B_C v = y$. From \EQ{wg0} and
\EQ{uIboundD}, we have
\begin{eqnarray*}
| v |^2_{\widetilde{A}}&=&| \u^{(I)}+ \u^{(II)}|^2_{H^1}  \le|
\u^{(I)} |^2_{H^1}+| \u^{(II)} |^2_{H^1} =\left(1+\frac{1}{\beta^2}\right)| \u^{(I)} |^2_{H^1}  +  \frac{C}{\beta^2 h^3}\left<g_{p_\Gamma},g_{p_\Gamma}\right>\\
& = & \left\{\begin{array}{cc}
\displaystyle{\left( 1 + \frac{1}{\beta^2} \right) | {\bf u}_{\Delta}^{(I)}
|^2_{A_{\Delta\Delta}} + \frac{C}{\beta^2
  h^3}\left<g_{p_\Gamma},g_{p_\Gamma}\right>}, &\mbox{for lumped
  preconditioner,} \\[2.4ex]
\displaystyle{\left( 1 + \frac{1}{\beta^2} \right) | {\bf u}_{\Delta}^{(I)}
|^2_{H_\Delta} + \frac{C}{\beta^2
  h^3}\left<g_{p_\Gamma},g_{p_\Gamma}\right>}, &\mbox{for Dirichlet
  preconditioner},\end{array}\right.
\end{eqnarray*}
where we used \EQ{Iprop} in the last equality.

On the other hand, we have from \EQ{precondlumped}--\EQ{subdomainwise}
\begin{eqnarray*}
\left< M^{-1}y,y \right>&=&
\frac{\myalpha}{h^3}\left<g_{p_\Gamma}, g_{p_\Gamma}\right>+ g_\lambda^T
M^{-1}_{\lambda} g_\lambda\\
& = & \left\{ \begin{array}{cc}
\displaystyle{\frac{\myalpha}{h^3}\left<g_{p_\Gamma}, g_{p_\Gamma}\right>+ g_\lambda^T B_{\Delta, D} A_{\Delta\Delta} B_{\Delta, D}^T g_\lambda}, & \mbox{for lumped preconditioner,}\\[2.4ex]
\displaystyle{\frac{\myalpha}{h^3}\left<g_{p_\Gamma}, g_{p_\Gamma}\right>+ g_\lambda^T
B_{\Delta, D} H_{\Delta} B_{\Delta, D}^T g_\lambda}, & \mbox{for Dirichlet preconditioner,}
\end{array}\right.
\\[2.4ex]
& = & \left\{ \begin{array}{cc}
\displaystyle{\frac{\myalpha}{h^3}\left<g_{p_\Gamma},g_{p_\Gamma}\right> + | {\bf u}_{\Delta}^{(I)} |^2_{A_{\Delta\Delta}}},  & \mbox{for lumped preconditioner,}\\[2.4ex]
\displaystyle{\frac{\myalpha}{h^3}\left<g_{p_\Gamma},g_{p_\Gamma}\right> + | {\bf u}_{\Delta}^{(I)} |^2_{H_\Delta}},  & \mbox{for Dirichlet preconditioner.}
\end{array}\right.
\end{eqnarray*}

It is not difficult to  see that $\left <v, v\right>_\Atilde \neq
0$. Otherwise, all the velocity components of $v$ would be zero,
cf. \EQ{wg0}, and then $B_C v$ would be zero, which conflicts with that $B_C v = y$ and $y$ is nonzero. \qquad \endproof

The proofs of the following two lemmas can be found at \cite[Lemmas 6.6 and 6.3]{li12}.
\begin{mylemma}
\label{lemma:BcW0}
For any $v = \left( {\bf w}_I, ~p_I, ~{\bf w}_{\Delta}, ~{\bf w}_{\Pi} \right) \in \Vtilde_0$, $B_C v \in R_G$.
\end{mylemma}

\begin{mylemma}
\label{lemma:m1R} For any $x \in R_{M^{-1} G}$,
\[
\left<Mx, x \right> = \max_{y \in R_G, y \neq 0} \frac{\left<y, x\right>^2}{\left<M^{-1}y,y\right>}.
\]
\end{mylemma}

The condition number bound of the preconditioned operator $M^{-1} G$
is given in the following theorem.

\begin{mytheorem}
\label{theorem:tcond} There exist positive constants $c$ and $C$, such that for all  $x \in R_{M^{-1} G}$,
\[
\min\left\{ 1, \myalpha\right\}  \frac{c \beta^2}{(1+\beta^2)} \left<Mx,x \right>\leq \left< G x,x \right> \leq C\left(\myalpha+ \Phi(H/h) \right)\left<
Mx,x \right>.
\]
\end{mytheorem}

\beginproof We only need to prove the above inequalities for any nonzero $x \in R_{M^{-1} G}$. We know from Lemma \LA{CG} that
\[
0 \neq \left< Gx,x\right> = x^T B_C \Atilde^{-1} B_C^Tx= x^T B_C \Atilde^{-1} \Atilde \Atilde^{-1} B_C^Tx =
\left< \Atilde^{-1} B_C^Tx, \Atilde^{-1} B_C^Tx\right>_{\Atilde}.
\]
Therefore $\Atilde^{-1} B_C^Tx \neq 0$. Also note that $\Atilde^{-1} B_C^Tx \in \Vtilde_{0}$ and $\left< \cdot , \cdot \right>_{\Atilde}$ defines a semi-inner product on $\Vtilde_{0}$, cf \EQ{wg0}, and then we have
\be \label{equation:Fnorm}
\left< Gx,x\right> = \max_{v \in \Vtilde_{0}, \left<v, v \right>_{\widetilde{A}} \neq 0} \frac{\left<v, \Atilde^{-1} B_C^Tx \right>^2_{\widetilde{A}}}{\left<v, v \right>_{\widetilde{A}}} =\max_{v \in \Vtilde_0, \left<v, v \right>_{\widetilde{A}} \neq 0}
\frac{\left<B_Cv,x\right>^2}{\left<\Atilde v,v\right>}.
\ee

{\it Lower bound:} From Lemma \ref{lemma:lower}, we know that for any nonzero $y \in R_G$, there exits
$w \in \Vtilde_{0}$, such that $B_C w = y$, $\left<w, w \right>_{\widetilde{A}} \neq 0$, $\left <\Atilde w, w\right> \le \max\left\{ 1,\frac{1}{\myalpha} \right\} \frac{C(1+\beta^2)}{\beta^2}  \left< M^{-1}y, y \right>$. Then from  \EQ{Fnorm}, we have
$$
\left< Gx,x\right> \ge\frac{\left<B_C w,x\right>^2}{\left<\Atilde w,w\right>}
\ge c \frac{\beta^2}{\max\left\{ 1,\frac{1}{\myalpha} \right\}(1+\beta^2)} \frac{\left<y,x\right>^2}{\left<M^{-1}y,y\right>}.
$$
Since $y$ is arbitrary, using Lemma \ref{lemma:m1R}, we have
$$\left< Gx,x\right> \ge c  \frac{\beta^2}{\max\left\{ 1,\frac{1}{\myalpha} \right\}(1+\beta^2)}  \max_{y \in R_G, y \neq 0} \frac{\left<y,x\right>^2}{\left<M^{-1}y,y\right>} = \min\left\{1,\myalpha\right\} \frac{c \beta^2}{(1+\beta^2)} \left<Mx,x \right>.$$

{\it Upper bound:}
From \EQ{Fnorm} and the fact that $\left< Gx,x\right> \neq 0$, we have
\[
\left< Gx,x\right> = \max_{v \in \Vtilde_0, \left<v, v \right>_{\widetilde{A}} \neq 0}
\frac{\left<B_Cv,x\right>^2}{\left<\Atilde v,v\right>} = \max_{v \in \Vtilde_0, \left<v, v \right>_{\widetilde{A}} \neq 0, B_Cv \neq 0}
\frac{\left<B_Cv,x\right>^2}{\left<\Atilde v,v\right>},
\]
where the maximum only needs to be considered among $v$ also satisfying $B_Cv \neq 0$. Then using
Lemmas \ref{lemma:upper}, \ref{lemma:BcW0}, and
\ref{lemma:m1R}, we have
\begin{eqnarray*}
\left< Gx,x\right> & \le & C (\myalpha+\Phi(H,h))\max_{v\in \Vtilde_{0}, \left<v, v \right>_{\widetilde{A}} \neq 0, B_Cv \neq 0} \frac{\left<B_C v,x\right>^2}{\left<M^{-1} B_C v,B_C v\right>} \\
& \le & C  (\myalpha+\Phi(H,h)) \max_{y\in R_G, y \neq 0} \frac{\left<y,x\right>^2}{\left<M^{-1}y,y\right>} = C  (\myalpha+\Phi(H,h)) \left<Mx,x\right>. \qquad \Box
\end{eqnarray*}

\begin{myremark}
\label{remark:alpha}
{\rm
We can see from Theorem \ref{theorem:tcond} that, for $\alpha \ge 1$,
the condition number bound of $M^{-1}G$  is proportional to
$\alpha+\Phi(H,h)$, and we should take smaller $\alpha$ to achieve faster convergence.
When $\alpha \leq 1$, the condition number bound is proportional to
$1+\frac{\Phi(H,h)}{\alpha}$ and we should take larger $\alpha$.  This explains why the value of $\alpha$ in \EQ{lumped} and \EQ{Dirichlet} is typically taken as~$1$. We introduce $\alpha$ in the preconditioner just for the convenience to demonstrate the convergence rates of the proposed algorithm in the following section.
}
\end{myremark}

\section{Numerical experiments}
\label{section:numerics}
We illustrate the convergence rate of the proposed algorithm by solving the incompressible Stokes problem \EQ{Stokes} in both two and three dimensions, on $\Omega = [0,1]^2$ and $\Omega=[0,1]^3$, respectively. Zero Dirichlet boundary condition is used. The right-hand side $\vvec{f}$ is chosen such that the exact solution is
$$\u=\left[\begin{array}{c}
\sin^3(\pi x)\sin^2(\pi y)\cos(\pi y)\\[0.8ex]
-\sin^2(\pi x)\sin^3(\pi y)\cos(\pi x)
\end{array}\right] , \quad
p=x^2-y^2,
$$
for two dimensions, and for three dimensions
$$\u=\left[\begin{array}{c}
\sin^2(\pi x)\left(\sin(2\pi y)\sin(\pi z)-\sin(\pi y)\sin(2\pi z)\right)\\[0.8ex]
\sin^2(\pi y)\left(\sin(2\pi z)\sin(\pi x)-\sin(\pi z)\sin(2\pi x)\right)\\[0.8ex]
\sin^2(\pi z)\left(\sin(2\pi x)\sin(\pi y)-\sin(\pi x)\sin(2\pi y)\right)
\end{array}\right], \quad
p=xyz-\frac{1}{8}. 
$$

The $Q_2$-$Q_1$ Taylor-Hood mixed finite element with continuous pressures is used; its inf-sup stability can be found at~\cite{ber, sou}. In two dimensions, the velocity space contains piecewise biquadratic functions and the pressure space contains piecewise bilinear functions; in three dimensions, piecewise triquadratic functions for the velocity and piecewise trilinear functions for the pressure.

The preconditioned conjugate gradient method is used to solve \EQ{Mspd}; the iteration is stopped when the $L^2-$norm of the residual is reduced by a factor of $10^{-6}$.

The following tables list the minimum and maximum eigenvalues of the
iteration matrix $M^{-1} G$, and the iteration counts for using both lumped
and Dirichlet preconditioners, respectively, for different cases. Here the extreme eigenvalues of $M^{-1} G$ are estimated by using the tridiagonal Lanczos matrix generated in the iteration.

\begin{table}[t]
\caption{\label{table:M1}
Performance of solving two-dimensional problem on $[0,1]^2$, $\alpha=1$ in \EQ{lumped} and \EQ{Dirichlet}.}
\centering
\begin{tabular}{cccccccccc}
\hline
&&&\multicolumn{3}{c}{lumped} & & \multicolumn{3}{c}{Dirichlet} \\
\cline{4-6} \cline{8-10}
\quad $H/h$  \quad & \#sub   &  & $\lambda_{min}$
& $\lambda_{max}$ & iteration  & & $\lambda_{min}$
& $\lambda_{max}$ & iteration \\ [.8ex]
\hline
$8$ &  $4  \times  4$ &&0.3066&32.28&31&&0.2983&4.40&18\\[.8ex]
    &  $8  \times  8$ &&0.3067&37.25&46&&0.2859&5.03&24\\[.8ex]
    & $16  \times  16$&&0.3068&38.42&51&&0.2556&5.28&25\\[.8ex]
    & $24  \times  24$&&0.3069&38.62&51&&0.2397&5.33&25\\[.8ex]
    & $32  \times  32$&&0.3070&38.68&51&&0.2304&5.36&25\\
\hline
& & & \\
\quad \#sub \quad & \quad $H/h$ \quad  & & $\lambda_{min}$
&  $\lambda_{max}$  & iteration & & $\lambda_{min}$
&  $\lambda_{max}$  & iteration  \\[.8ex]
\hline
$8\times 8$  & $4  $ &&0.3024&15.91&34 &&0.2706&4.15&21\\[.8ex]
             & $8  $ &&0.3067&37.25&46 &&0.2859&5.03&24\\[.8ex]
             & $16 $ &&0.3069&85.32&62 &&0.2966&6.04& 25\\[.8ex]
             & $24 $ &&0.3073&137.49&73&&0.3028&6.69&26\\[.8ex]
             & $32 $ &&0.3075&192.32&83&&0.3070&7.19  &27\\
\hline
\end{tabular}
\end{table}

Table \ref{table:M1} shows the performance for solving the two-dimensional problem. The coarse level velocity space in the algorithm is spanned by the subdomain vertex nodal basis functions corresponding to each velocity component. We take $\myalpha=1$ in both the lumped and the Dirichlet preconditioners \EQ{lumped} and \EQ{Dirichlet}.  We can see from Table \ref{table:M1} that the minimum eigenvalue is independent of the mesh
size for both preconditioners. The maximum eigenvalue is independent of the
number of subdomains for fixed $H/h$; for fixed number of subdomains,
it depends on $H/h$ in the order of $(H/h) (1 + \log{(H/h)})$ for the
lumped preconditioner, and in the order of $(1 + \log{(H/h)})^2$ for
the Dirichlet  preconditioner.

Tables \ref{table:M2} and \ref{table:M3} are for solving the
three-dimensional problem. The coarse level velocity space is spanned
by the subdomain vertex nodal basis functions and subdomain edge-cutoff functions corresponding to  each velocity component. This coarse space is the same as for solving scalar elliptic
problems in \cite[Algorithm 6.25]{Toselli:2004:DDM} corresponding to
each velocity component. In Table \ref{table:M2}, $\myalpha=1$; in Table \ref{table:M3}, $\myalpha=1/2$.

In Table \ref{table:M2}, the minimum eigenvalue is independent of the mesh
size for both preconditioners. The maximum eigenvalue is independent of the
number of subdomains for fixed $H/h$; for fixed number of subdomains,
it depends on $H/h$, but not in the order of $(H/h) (1 + \log{(H/h)})$ for the
lumped preconditioner, nor  $\left(1 + \log{(H/h)}\right)^2$ for the
Dirichlet  preconditioner, as $\Phi(H/h)$ does.  Moreover, the
convergence rate of the algorithm using the Dirichlet preconditioner
is only slightly better than using the lumped preconditioner. The reason is that the upper eigenvalue bound in Theorem \ref{theorem:tcond} depends on two terms $\alpha$ and $\Phi(H/h)$, and in this case $\alpha = 1$ dominates when $H/h$ is small. Therefore, even though using the Dirichlet preconditioner can reduce $\Phi(H/h)$ compared with using the lumped preconditioner, this improvement on the upper eigenvalue bound can not show up in Table \ref{table:M2}. What shows in Table \ref{table:M2} for $\lambda_{max}$ is essentially its dependence on $\alpha$. Only for larger $H/h$, e.g., for $H/h=6$ and $H/h=8$ in Table \ref{table:M2}, the improvement on the upper eigenvalue bound by using the Dirichlet preconditioner becomes visible.
%Due to the memory limitation, we could not run the experiments with large $H/h$ when we change the number of subdomains.

To experiment the case when $\alpha$ is less dominant in the upper eigenvalue bound, we take $\alpha = 1/2$ in Table \ref{table:M3}. Consistent with Theorem \ref{theorem:tcond}, the lower eigenvalue bounds in Table \ref{table:M3} become half of those in Table  \ref{table:M2} and they are also independent of the mesh size.  The upper eigenvalue bounds exhibit the pattern of $\Phi(H/h)$ for both preconditioners. They are independent of the
number of subdomains for fixed $H/h$; for fixed number of subdomains,
they depend on $H/h$ in the order of $(H/h) (1 + \log{(H/h)})$ for the
lumped preconditioner, and in the order of $(1 + \log{(H/h)})^2$ for the
Dirichlet  preconditioner.

\begin{table}[t]
\caption{\label{table:M2}
Performance of solving three-dimensional problem on $[0,1]^3$, $\alpha=1$ in \EQ{lumped} and \EQ{Dirichlet}.}
\centering
\begin{tabular}{cccccccccc}
\hline
&&&\multicolumn{3}{c}{lumped} & & \multicolumn{3}{c}{Dirichlet} \\
\cline{4-6} \cline{8-10}
\quad $H/h$ \quad & \#sub   & & $\lambda_{min}$
& $\lambda_{max}$ & iteration  & & $\lambda_{min}$
& $\lambda_{max}$ & iteration \\[.8ex]
\hline
$4$ &  $3 \times  3 \times  3$&&0.0776&9.13&56&&0.0776&8.97&56\\[.8ex]
    &  $4 \times  4 \times  4$&&0.0775&9.35&54&&0.0774&9.19&55\\[.8ex]
    &  $6 \times  6 \times  6$&&0.0773&9.41&58&&0.0773&9.23&59\\[.8ex]
    &  $8 \times  8 \times  8$&&0.0773&9.51&57&&0.0772&9.34&61\\
\hline
& & & \\
\quad \#sub \quad & \quad $H/h$ \quad  & & $\lambda_{min}$
&  $\lambda_{max}$  & iteration & & $\lambda_{min}$
&  $\lambda_{max}$  & iteration  \\[.8ex]
\hline
$3\times 3 \times 3$  & $3 $ & &0.0760&8.06&54 &&0.0760&7.96&54\\[.8ex]
                   & $4 $ & &0.0776&9.13&56 &&0.0776&8.97&56 \\[.8ex]
                   & $6 $ & &0.0780&11.88&53&&0.0780&9.35& 55\\[.8ex]
                   & $8 $ & &0.0780&16.64&57&&0.0780&9.44&55\\
\hline
\end{tabular}
\end{table}

\begin{table}[t]
\caption{\label{table:M3} Performance of solving three-dimensional problem on $[0,1]^3$, $\alpha=1/2$ in \EQ{lumped} and \EQ{Dirichlet}.}
\centering
\begin{tabular}{cccccccccc}
\hline
&&&\multicolumn{3}{c}{lumped} & & \multicolumn{3}{c}{Dirichlet} \\
\cline{4-6} \cline{8-10}
\quad $H/h$ \quad & \#sub   &  & $\lambda_{min}$
& $\lambda_{max}$ & iteration  & &  $\lambda_{min}$
& $\lambda_{max}$ & iteration \\[.8ex]
\hline
$4$ & $3 \times  3 \times  3$ &&0.0395&7.20&59&&0.0395&4.89&54\\[.8ex]
    & $4 \times  4 \times  4$ &&0.0394&8.15&66&&0.0394&5.01&53\\[.8ex]
    & $6 \times  6 \times  6$ &&0.0393&8.85&70&&0.0393&5.03&55\\[.8ex]
    & $8 \times  8 \times  8$ &&0.0393&9.09&72&&0.0393&5.09&56\\
\hline
& & & \\
\quad \#sub \quad & \quad $H/h$ \quad  & & $\lambda_{min}$
&  $\lambda_{max}$  & iteration & &  $\lambda_{min}$
&  $\lambda_{max}$  & iteration  \\[.8ex]
\hline
$3\times 3 \times  3$  & $3 $ &&0.0387&5.15&55 &&0.0387&4.35&53\\[.8ex]
                       & $4 $ &&0.0395&7.20&57 &&0.0395&4.89&54\\[.8ex]
                       & $6 $ &&0.0397&11.70&63&&0.0397&5.11& 52\\[.8ex]
                       & $8 $ &&0.0397&16.52&73&&0.0397&5.17&52\\
\hline
\end{tabular}
\end{table}

\end{document}